\documentclass[11pt]{article}
\usepackage{srcltx}
\usepackage{times}
\usepackage{amsmath,amssymb,amsfonts}
\newcounter{theorem}
\newcounter{lemma}

\def\@cite#1#2{[{\bf #1\if@tempswa , #2\fi}]}
\def\qed{\ifhmode\unskip\nobreak\fi\ifmmode\ifinner\else\hskip.5em\fi\fi
 \hbox{\hskip.5em$\square$\hskip.1em}}
\begin{document}
\title{\sc Hardy-type results on the average of the lattice point error term over long intervals}
\author{Burton Randol}
\date{}
\maketitle
\bibliographystyle{plain}

\begin{abstract}\begin{sloppypar}Suppose $D$ is a suitably admissible compact subset of $\mathbb{R}^k$ having a smooth boundary with possible zones of zero curvature. Let \mbox{$R(T,\theta,x)= N(T,\theta,x) - T^{k}\mathrm{vol}(D)$,} where $N(T,\theta,x)$ is the number of integral lattice points contained in an $x$-translation of $T\theta(D)$, with $T  >0$ a dilation parameter and $\theta \in SO(k)$. Then $R(T,\theta,x)$ can be regarded as a function with parameter $T$ on the space $E_{*}^{+}(k)$, where $E_{*}^{+}(k)$ is the quotient of the direct Euclidean group by the subgroup of integral translations, and $E_{*}^{+}(k)$ has a normalized invariant measure which is the product of normalized measures on $SO(k)$ and the $k$-torus. We derive an integral estimate, valid for almost all $(\theta,x) \in E_{*}^{+}(k)$, one consequence of which in two dimensions is that for almost all $(\theta,x) \in E_{*}^{+}(2)$, a counterpart of the Hardy circle estimate \mbox{$(1/T)\int_{1}^{T} |R(t,\theta,x)\,dt| \ll T^{\frac{1}{4} +\epsilon}\;$}is valid with an improved estimate. We conclude with an account of hyperbolic versions for which, drawing on previous work of Hill and Parnovski \cite{hill-parnovski}, we give counterparts in all dimensions, for both the compact and non-compact finite volume cases.\end{sloppypar}\end{abstract}

\bigskip

Define $N(T)$ to be the number of integral lattice points in $\mathbb{R}^2$ which are at distance $\leq T$ from the origin. Then the study of the so-called remainder term $N(T)-\pi T^2$, the difference of $N(T)$ from the area of a disk of radius $T$, has been pursued for many decades. Various ways of measuring the size, for large $T$, of $R(T)$ have been investigated, for example, asymptotic estimates for $R(T)$ itself, as well as estimates on various averages, for instance, \[\frac{1}{T}\int_1^T |R(t)|\,dt\,.\] There are also corresponding results when the center is allowed to shift, and become a new parameter in the definition of $R$ (cf.\ \cite{bleher2}, \cite{bleher1}, \cite{kendall},). We note that it was usual in the older literature to replace $T$ by $\sqrt{T}$ in the definitions of $N(T)$ and $R(T)$, since then $N(T)$ counts the number of ordered pairs $(n_1,n_2)$ in $\mathbb{Z}\times\mathbb{Z}$ with $n_1^2 + n_2^2 \leq T$. The definition with $T$ in place of $\sqrt{T}$ is generally more usual if one regards the question geometrically, as arising from dilations of a fixed body $D$ (in this case the closed unit disk around the origin).

Expressed with the $\sqrt{T}$ scaling, an important result of the averaging type was obtained by Hardy in \cite{hardy1}, in which he showed that for any $\epsilon > 0$,\[\frac{1}{T} \int_1^T |R(t)|\,dt \ll T^{\frac{1}{4} + \epsilon}\,,\] a result which was subsequently refined by Cram\'{e}r, in \cite{cramer}, with later improvements over time by numerous authors.

Hardy, at the beginning of \cite{hardy1}, wrote that it is ``not unlikely'' that for any $\epsilon >0$, \[R(T) \ll T^{\frac{1}{4} + \epsilon}\,.\] This estimate has not yet been confirmed or refuted, although Hardy's average result is suggestive of its possible truth.

The originating question has led to numerous fruitful areas of investigation, e.g., the asymptotic study of the lattice point count in dilations of a given general domain in Euclidean space $\mathbb{R}^k$, or the asymptotics of the count, within a dilating ball in hyperbolic $k$-space $\mathbb{H}^k$, of elements of the orbit$\{\gamma(x)\}$ ($\gamma \in \Gamma$) of a point under the action of a discontinuous group $\Gamma$. We note that the Euclidean question can be directly related to counting the number of eigenvalues $\leq T$ for differential operators on the integral $k$-torus, i.e., the quotient of $\mathbb{R}^k$ by $\mathbb{Z}^k$. For example, if $P(x_1,\ldots,x_k)$ is a real polynomial which is homogeneous of degree $h$ and positive except at the origin, the operator $P(\frac{\partial}{\partial x_1}\ldots ,\frac{\partial}{\partial x_k})$ has eigenvalues $P((2\pi i n_1,\ldots ,2\pi i n_k) = (2\pi i)^h P((n_1,\ldots ,n_k)$ for $(n_1,\ldots ,n_k) \in \mathbb{Z}^k$, and so up to a real factor (by positivity of $P$, $h$ must be even), one is counting the number of lattice points in dilates by $T^{1/h}$ of the domain defined by $P(x_1,\ldots,x_k) \leq 1$, with the convention that we count eigenvalues of the negative of the operator when $4$ does not divide $h$. 

There are counterparts for general Riemannian manifolds of the eigenvalue asymptotic question. To take an example, suppose $S$ is a compact Riemann surface of genus greater than $1$ and constant curvature $-1$, and $0=\lambda_0 < \lambda_1 \leq \lambda_2 \ldots$ is the sequence of eigenvalues for the problem $\Delta f + \lambda f =0$ on $S$. By the Weyl asymptotic law, if $N(T)$ is the corresponding eigenvalue count in $[0,T]$, $N(T) = (A/4\pi)T + R(T)$, where $A$ is the area of the surface, and the remainder $R(T)$ is $o(T)$. This general estimate on $R(T)$ can be considerably improved (cf.\ \cite{hejhal1}, \cite{randol78}), and, as in the case of the circle problem, examination of the analysis makes it tempting to believe that $R(T) \ll T^{\frac{1}{4} + \epsilon}\,$ for any $\epsilon > 0$. It has been shown in arithmetic cases (cf.\ \cite{hejhal1}) that the $\frac14$ in this estimate cannot be reduced. I have not been able to establish or refute the estimate, but in \cite{randol81} was able to establish its almost everywhere validity for a closely related version of the eigenvalue remainder problem, in which the eigenfunctions figure in the count.

Intrigued by this result, LaPointe, Polterovich, and Safarov in \cite{polterovich1} subsequently obtained several interesting results of this type in a general Riemannian context.

After this brief introduction to the general background and context for the results of this paper, we return to the theme of asymptotic estimates for lattice point counts in dilating bodies in $\mathbb{R}^k$, in particular, to considerations in which rotations as well as shifts are taken into account. The inclusion of rotations becomes necessary in the presence of zones of zero curvature on $\partial D$, because the asymptotics in this case can be exquisitely sensitive to the rotational orientation of $D$, and lattice point asymptotics in the zero curvature case can sharply differ from those in the positive curvature case (cf.\cite{randol66a}, \cite{randol66b}). An important example where vanishing curvature can play a role arises in the case of positive homogeneous forms that are not quadratic, for example, $x_1^{2m} + \cdots x_k^{2m}$, with $m > 1$. 

Accordingly, we will now examine the lattice point count corresponding to dilations of a rotation $\theta(D)$ of the body $D$, with $\theta \in SO(k)$.

Since the integral lattice point count over $T\theta(D)$ is clearly unaffected by integral translations, it can be naturally regarded as a parametrized function $N(T,\theta,x)$ on the integral torus $\mathbb{T}^k$, with $x \in \mathbb{T}^k$.

The Poisson summation formula is true in this context as a statement in $L^2(\mathbb{T}^k)$, and shows that in $L^2(\mathbb{T}^k)$, \[N(T,\theta,x) = \sum_n \hat{\chi}_T(\theta^{-1}(n))\,e^{2\pi (n,x)}\,,\] where $\hat{\chi}_T$ is the Fourier transform of the indicator function $\chi_{T}$ of $TD$, and the summation is over $\mathbb{Z}^k$. Since $\hat{\chi}_T(\theta^{-1}(n)) = T^k\hat{\chi}(T\theta^{-1}(n))$, with $\hat{\chi} = \hat{\chi}_{1}$, and $\hat{\chi}(0) = \mathrm{vol}(D)$, we thus find that in  $L^2(\mathbb{T}^k)$, the remainder term, $R(T,\theta,x) = N(T,\theta,x) - T^k\mathrm{vol}(D)$ equals
\[T^k\sideset{}{'}\sum \hat{\chi}(T\theta^{-1}(n))\,e^{2\pi (n,x)}\,,\] where the prime indicates that the origin is omitted from the sum.

At this point, we briefly review some results which establish comparison theorems between the decay of the Fourier transform of the indicator function of a compact domain $D$ in $\mathbb{R}^k$ and the decay of the Fourier transform of the indicator function of a ball. Of particular interest, if $\partial D$ is adequately smooth, are examples having zones of zero curvature on $\partial D$. We begin with an early pointwise estimate, stronger than required for the specific purposes of this paper, but applicable to some situations to which the later $L^2$ averaging estimates are not. It was motivated by questions arising from the study of representations by positive forms.

\medskip

{\sc Lemma A (\cite{randol69a}, \cite{randol69b})}.{\sl \ Suppose $D$ is a body in $\mathbb{R}^k$, satisfying

\begin{description}

\item[If $\mathbf{k=2}$:]  $\partial D$ is of class $C^{n+3}$ for some integer $n \geq 1$, and the Gaussian
curvature of $\partial D$ is nonzero at all points of $\partial D$, with the possible exception of a finite
set, at each point of which the tangent line has contact of order $\leq n$. 

\item[If $\mathbf{k\geq 3}$:] $D$ is convex and $\partial D$ is real analytic. 

\end{description}

Let $\hat{\chi}(r,\phi)$ be the Fourier transform, in polar coordinates of the indicator function of $D$. Then \[|\hat{\chi}(r,\phi)| \ll \Psi(\phi)r^{-(n+1)/2}\,,\]where $\Psi \in L^p(S^{k-1})$, for some (in principle computable) $p >2$. As a pointwise estimate about numbers, this is of course to be interpreted as valid for almost all $\phi \in S^{k-1}$.}

\medskip

In other words, the Fourier transform of a body in either of the above classes decays like that of a ball, up to multiplication by an $L^p$ function on the sphere, with $p >2$. Note that this result implies the weaker corollary that for some $p > 2$,\[\int_{S^{k-1}}|\hat{\chi}(r,\phi)|^p\,d\phi \ll r^{-(p/2)(k+1)}\,.\] 

\smallskip

{\sc Remark.} The real analyticity can be relaxed, as Svennson, following a suggestion of H\"{o}rmander, showed by extending the above result to the $C^{\infty}$ category, with an additional hypothesis about the absence of tangent lines having infinite contact order with $\partial D$  \cite{svensson}, but I do not know the current extent to which the convexity and contact requirements can be relaxed for this type of pointwise result (the possibility that convexity might be inessential was suggested in \cite{randol69b}). Both requirements can be eliminated for some $L^2$ averaging results (see below), for which a good recent reference and bibliography is \cite{iosevich}. 

\smallskip

As noted above, for results of $L^2$ averaging type, more general types of domains can be handled. We now describe such a result, valid for $C^{\infty}$ boundaries without a convexity or contact condition, which arose from a conjecture Arnold put forward in his seminar during the course of expositions of Lemma A and a later result of Colin de Verdi\`{e}re \cite{verdiere1}, and which was established in final form by Varchenko. The removal of the convexity requirement in the resulting averaging theorem unrestrictedly extends the applicability of the results of the present paper to all positive homogeneous forms.

\medskip

{\sc Lemma B (\cite{varchenko1}, \cite{varchenko3})}.{\sl \  Suppose $D$ is compact in $\mathbb{R}^k$, with $C^{\infty}$ boundary. Then \[\int_{S^{k-1}}|\hat{\chi}(r,\phi)|^2\,d\phi \ll r^{-(k+1)}\,.\] } 

\medskip

Later, Brandolini, Hofmann, and Iosevich, as part of an ongoing program to explore the maximum generality of such averaging results, established the above $L^2$ estimate for general bounded convex domains in $\mathbb{R}^k$, and for domains without the convexity restriction, but satisfying a $C^{3/2}$ smoothness hypothesis. 

\medskip

{\sc Lemma C (\cite{brandolini2})}.{\sl \ If $D$ is convex, or of class $C^{3/2}$ with no convexity hypothesis, then \[\int_{S^{k-1}}|\hat{\chi}(r,\phi)|^2\,d\phi \ll r^{-(k+1)}\,.\]}

\medskip

We are now ready to establish the following theorem, taking normalized Haar measure on $SO(k)$ and on $\mathbb{T}^k$, and the invariant product measure on the space of pairs $E_{*}^{+}(k) = (\theta,x)$, which is the quotient of the direct Euclidean group by the subgroup of integral translations.

{\sc Theorem 1}. {\sl  Suppose $f(t)> 0$ is differentiable and non-increasing on $[a,\infty)$, with $a>0$ and $\int_a^\infty f(t)/t\,dt < \infty$, and assume that $D$ is a body in $\mathbb{R}^k$ of a type described in the lemmas. Then for almost all $(\theta,x)$, \[\frac{1}{T}\int_a^T |R(t,\theta,x)|\,dt \ll T^{(k-1)/2} (f(T))^{-1/2}\,,\] where the implied constant in the inequality depends on the particular choice of $\theta$ and $x$.} 

{\sc Corollary}.  {\sl (Counterpart of Hardy's theorem). For almost all $(\theta,x)$ and for any $\epsilon >0$, \[\frac{1}{T}\int_1^T |R(t,\theta,x)|\,dt \ll T^{((k-1)/2) + \epsilon\,.}\]}(Take $f(t) = t^{-2\epsilon}$).

{\sc Remark.} If we scale by $T^{1/h}$ instead of by $T$, $T^{1/2}$ being the scaling used by Hardy in his treatment of the circle case, the estimate on the right becomes $T^{((k-1)/2h)+\epsilon}$, which for $k=2$ corresponds to the Hardy estimate $T^{\frac{1}{4} + \epsilon}$.

\medskip

It is obvious from the theorem that stronger results are true. For example,

\medskip

{\sc Corollary}. {\sl  For almost all $(\theta,x)$ and any $\epsilon > 0$, \[\frac{1}{T}\int_2^T |R(t,\theta,x)|\,dt \ll T^{(k-1)/2}\log^{\frac{1}{2} +\epsilon}T \,.\]}(take $f(t) = \log^{-(1+\epsilon)}t$).

\medskip

{\sc Proof of Theorem 1}. As we have seen, in $L^2(\mathbb{T}^k)$, \[R(T,\theta,x) = T^k\sideset{}{'}\sum \hat{\chi}(T\theta^{-1}(n))\,e^{2\pi (n,x)}\,,\]
so the result of integrating $|R(T,\theta,x)|^2$ over the torus is \[T^{2k} \sideset{}{'}\sum |\hat{\chi}(T\theta^{-1}(n))|^2\,.\]

Now invoking the dimensionally appropriate estimate and noting that $|T\theta^{-1}(n)| = T|n|$, we see that this quantity can be bounded by \[T^{k-1}\sideset{}{'}\sum |\Psi(\theta^{-1}(n_{\phi})|^2|n|^{-(k+1)}\,,\]where $\Psi \in L^2(S^{k-1})$, and $n_{\phi}$ is the spherical component of $n$.

Integrating over $SO(k)$, noting that the result is the same as integrating $|\Psi|^2$ over the measure-normalized sphere, and that the series in the resulting estimate converges absolutely, we find that \[ \int_{SO(k)}d\theta\int_{\mathbb{T}^k} |R(T,\theta,x)|^2\, dx \ll T^{k-1}\,.\]

We next remark that throughout what follows, the dependency on parameters of constants that will be associated with various estimates will not be explicitly mentioned, but should be clear from context. Bearing this in mind, note that the last estimate implies that \[\int_a^\infty t^{-k} f(t) \,dt\int_{SO(k)}d\theta\int_{\mathbb{T}^k} |R(t,\theta,x)|^2\,dx  \leq  C\,,\]for some $C >0$.

Now, as in \cite{randol81}, we note that by Fubini's theorem, the integral in the $t$ variable is finite for almost all $(\theta,x) \in E^+_*(k)$, which implies that for almost all \mbox{$(\theta,x) \in E^+_*(k)$,}\[\int_a^T t^{-k}f(t) |R(t,\theta,x)|^2\,dt \ll 1\,.\]

Define \[R_1(T,\theta,x) = \int_a^T |R(t,\theta,x)|^2 \,dt\,.\]

Then, integrating by parts, and noting that $R_1(a,\theta,x)= 0$,\[1 \gg\int_a^T t^{-k}f(t) |R(t,\theta,x)|^2\,dt =\] \[T^{-k}f(T) R_1(T,\theta,x) + \int_a^T (kt^{-k -1}f(t) - t^{-k}f'(t)) R_1(t,\theta,x)\,dt \,. \]

The second term in the last expression is non-negative since $f'(t) \leq 0$, which implies that

\[T^{-k}f(T) R_1(T,\theta,x) \ll 1\,,\]so \[ R_1(T,\theta,x) = \int_a^T |R(t,\theta,x)|^2 \,dt  \ll T^{k}(f(T))^{-1}\,.\]

\begin{sloppypar}This in turn implies by the Schwartz inequality, since $|R(t,\theta,x)| = \mbox{$1\times |R(t,\theta,x)|$}$, that \[\int_a^T |R(t,\theta,x)|\,dt \ll T^{1/2}T^{k/2}(f(T))^{-1/2}
 = T^{(k+1)/2}(f(T))^{-1/2}\,,\]so for almost all $(\theta,x)\in E^+_*(k)$,\[\frac{1}{T}\int_a^T |R(t,\theta,x)|\,dt \ll T^{(k-1)/2}(f(T))^{-1/2}\,.\qed\]\end{sloppypar}

{\sc Remark.} It can be shown by, for example, adapting the stationary phase argument used to derive an $\Omega$-estimate in \cite{randol66a}, that in the presence of zero curvature, the almost everywhere estimate of the theorem can fail for specific values of $(\theta,x) $.

\bigskip

We continue with a description of counterparts of these results in the hyperbolic case, for which the analytically natural dilating object is a ball (the Selberg theory, which furnishes the applicable method, deals with point-pair invariants, which are functions of distance). With this in mind, the goal is to first formulate a suitable hyperbolic analogue to the result in Kendall's foundational paper \cite{kendall} on variance of the lattice point count over shifted ovals. The formulation we will first consider is one in which $\Gamma$ is a torsion-free co-compact group acting on hyperbolic $k$-space $\mathbb{H}^k$, with the role of lattice points played by the $\Gamma$-orbit of a point $y$ in $\mathbb{H}^k$, the aim being to obtain a Hardy-type integral estimate for the lattice point count of $T$-dilates of a unit ball with center $x$, valid for almost all $(x,y) \in S \times S$, where $S$ is the compact hyperbolic manifold of constant curvature $-1$ and volume $V$ that is the quotient of $\mathbb{H}^k$ by $\Gamma$. Other questions of this type have been been considered by several authors, e.g., \cite{hill-parnovski}, \cite{huber56}, \cite{lax-phillips}, \cite{patterson}, \cite{phillips-rudnick}, \cite{randolchapter}, \cite{wolfe}. A discussion of a version in which $y$ is fixed but $x$ varies and the variance over the $x$ variable is studied can be found in \cite{hill-parnovski}, the latter formulation being more suited to the non-compact case, since in that case the lattice point count is not in $L^2(S \times S)$, as is noted in \cite{hill-parnovski} and \cite{wolfe}. We remark for later use that our general technique for deriving a.e.\ integral estimates is applicable to the results in \cite{hill-parnovski}. We also remark that the study of variance removes a barrier to the known accuracy of estimates in $T$ that are valid for specific points in $S \times S$, the analysis of which may not include possible contributions from small eigenvalues near $((k-1)/2)^2$ (cf.\ \cite{lax-phillips}, \cite{patterson}, \cite{phillips-rudnick}).

As the last sentence suggests, the results of the analysis depend very much on whether or not $S$ has small eigenvalues, i.e., on whether or not there are eigenvalues in $(0,((k-1)/2)^2)$ for the problem $\Delta f + \lambda f =0$ on $S$. It is known that any compact hyperbolic manifold with non-vanishing first Betti number has finite covers for which they are present, as close to $0$ as desired, and in numbers that can be made arbitrarily large by passage to sufficiently deep covers (cf.\ \cite{buser4g}, \cite{randol74}, \cite{randolchapter}). All compact hyperbolic surfaces satisfy the Betti number condition, and Millson \cite{millson} has demonstrated the existence of examples with non-vanishing first Betti number in all dimensions $\geq 3$. In particular, all hyperbolic compact Riemann surfaces have finite covers with small eigenvalues, the number of which is, however, always bounded by $2g-3$, where $g$ is the genus of the surface \cite{otal-rosas}, although it is possible that such eigenvalues do not occur for certain arithmetic surfaces. For instance, it is known in some cases that none are present in intervals of the form $(a,\frac{1}{4})$, with $a < \frac{1}{4}$ dependent on the surfaces being studied (cf., for example, \cite{luo-rudnick-sarnak}), but a very general result of this type that completely excludes $(0,\frac{1}{4})$ has not been established.

We will begin with a quick analysis of the illustrative example of $\mathbb{H}^3$, since the general ideas and techniques are exceptionally easy to follow in that case. In particular, the required transforms take a very simple form there and are explicitly given on \mbox{p.\ 280} of \cite{randolchapter}. We will afterwards discuss the general case, as well as the previously mentioned variant of the problem that is applicable to the non-compact finite volume case. We note that there is an inconsequential misprint on p.\ 280 of \cite{randolchapter}, on line 13 from the top, where what should be $O(r^{-2}e(t))$ is mistakenly printed as $O(r^{-3}e(t))$.

Suppose now that $\Gamma$ is as above, and that $N(T,x,y)$ is the number of points of the form $\gamma y$ ($\gamma \in \Gamma$, $y \in \mathbb{H}^k$) contained in the closed hyperbolic ball of radius $T$ centered at $x \in \mathbb{H}^k$. This can clearly be regarded as a function on $S \times S$, and in order to reproduce our previous argument we need an expression for the integral over $S \times S$ of a suitable quantity related to $N(T,x,y)$. The Selberg pre-trace formula is well adapted to this situation, and as mentioned, the necessary analysis is carried out in \cite{randolchapter}. In general, in $\mathbb{H}^k$, the pre-trace formula tells us that in $L^2(S\times S)$,\[N(T,x,y) = \sum_{n=0}^{\infty} h_T(r_n)\varphi_n(x)\varphi_n(y) \,,\] where $h_T(r)$ denotes the Selberg transform of the point-pair invariant corresponding to the indicator function of a ball of radius $T$. Here, as usual for general $k$, the $r_n$'s are related to the eigenvalues by $r_n = (\lambda_n -((k-1)/2))^{2})^{1/2}$ and the $\varphi_n$'s are corresponding orthonormal eigenfunctions of the Laplace operator on $S$. This correspondence produces two $r_n$'s for each eigenvalue, counting the $r_n$ at $0$ with multiplicity $2$ if it is present. In order to ensure a $1-1$ correspondence between eigenvalues and $r_n$'s,  we will only count those $r_n$'s that lie on the union of $(0,\infty)$ with the closed segment of the imaginary axis from the origin to $((k-1)/2)i$, and count the $r_n$ at $0$ once if it is present. The $r_n$'s on this segment of the imaginary axis correspond to small eigenvalues, augmented by the eigenvalue at $((k-1)/2)^2$ if present, which corresponds to $r_n = 0$, together with the eigenvalue $\lambda_0 =0$, which corresponds to $r_0 = ((k-1)/2)i$. For purposes of this paper, we will call all such eigenvalues special eigenvalues.

The above series for $N(T,x,y)$ in \mbox{$L^2(S\times S)$} contains a finite number of explicit terms coming from Selberg transforms corresponding to special eigenvalues, and infinitely many coming from the remaining eigenvalues corresponding to $r_n$'s on the real axis minus the origin. I.e., in $L^2(S\times S)$, we can rewrite our previous identity as\[N(T,x,y) -  \sum_{n=0}^{N} h_T(r_n)\varphi_n(x)\varphi_n(y) = \sum_{n=N+1}^{\infty} h_T(r_n)\varphi_n(x)\varphi_n(y) \,,\] where the sum on the left side is the contribution to the pre-trace formula from the special eigenvalues. Note that $\varphi_0(\cdot) =V^{-1/2}$, \ so \[h_T(r_0)\varphi_0(x)\varphi_0(y) = V^{-1}h_T(r_0)\,,\] and since it can be shown that $h_T(r_0) = B_k(T)$, where $B_k(T)$ is the volume of a ball of radius $T$ in $\mathbb{H}^k$, we can further rewrite the identity as\[N(T,x,y) - B(T)/V - \sum_{n=1}^{N} h_T(r_n)\varphi_n(x)\varphi_n(y) = \sum_{n=N+1}^{\infty} h_T(r_n)\varphi_n(x)\varphi_n(y) \,.\]

As mentioned, the description of the Selberg transform $h_T(r)$ is particularly simple for $\mathbb{H}^3$, and is as follows:

\begin{enumerate}

\item For all $r$ except $0$ and $i$, \[h_T(r) = (4\pi /r (1+r^2))(\cosh T \sin rT - r\sinh T \cos rT)\,.\]

\item $h_T(i) = \pi\sinh 2T - 2\pi T = B_3(T)$. 

\item $h_T(0) = 4\pi(T\cosh T -\sinh T)$.

\end{enumerate}

The following lemma is now a consequence of these facts, and can be inferred from either Parseval's theorem or Weyl's law for the $r_n$'s, specifically, in the case of the latter, the consequence that \[\sum_{n=N+1}^\infty r_n^{-4} < \infty.\]

\bigskip

{\sc Lemma.} {\sl Define, in $L^2(S\times S)$, \[\tilde{N}(T,x,y) = N(T,x,y) - B_3(T)/V\]\[ - (h_T(r_1)\varphi_1(x)\varphi_1(y) +\cdots + h_T(r_N)\varphi_N(x)\varphi_N(y))\]\[ = \sum_{n=N+1}^{\infty} h_T(r_n)\varphi_n(x)\varphi_n(y)\,.\] Then\[\int_{S\times S}|\tilde{N}(T,x,y)|^2 dxdy = \sum_{n=N+1}^\infty |h_T(r_n)|^2 \ll e^{2T}\,,\]}

\medskip

This estimate is employed exactly as in the treatment of the Euclidean case. Namely, we integrate it from $a$ to $\infty$ over the dilation parameter, using a measure suitably weighted to produce integrability, and then invoke a variant of the Fubini-Tonelli theorem to ultimately derive the local property of integrability from $a$ to $\infty$ in the dilation parameter, valid for almost all points in the base space \mbox{$S \times S$,} and  consequently boundedness, for almost any point in the base space, of the integrals from $a$ to $T$, where the bound depends on the point in the base space. This leads to an $L^2$ Hardy-type result, which is then converted to an $L^1$ result via the Schwartz inequality. This technique, which was used for a similar purpose in \cite{randol81}, and has undoubtedly been employed multiple times in various guises, is probably applicable to a wide variety of genericity questions in probabilistic number theory and other areas.

\medskip

{\sc Remark.} The last estimate can be supplemented with an estimate from below. In more detail, since for $r \in R^1 - \{0\}$,\[h_T(r) = (4\pi /r (1+r^2))(\cosh T \sin rT - r\sinh T \cos rT)\,,\]it follows that\[\int_{S\times S}|\tilde{N}(T,x,y)e^{-T}|^2 dxdy = \sum_{n=N+1}^\infty |e^{-T}h_T(r_n)|^2\]can be written in the form $\Psi(T) + \Phi(t)$, where $\Psi(T)$, a uniform limit of generalized trigonometric polynomials, is a non-vanishing Bohr almost periodic function, and $\Phi(T) \ll e^{-2T}$. In particular therefore, \[\int_{S\times S}|\tilde{N}(T,x,y)e^{-T}|^2 dxdy = \Omega(1)\,,\]or 
\[\int_{S\times S}|\tilde{N}(T,x,y)|^2 dxdy= \Omega(e^{2T})\,.\]

\medskip

This said, we return to the statement of our conclusion in the 3-dimensional hyperbolic case, and as noted, since the derivation is methodologically identical to that of the already discussed counterparts in the Euclidean case, we will simply state the result: 

\medskip

{\sc Theorem 2.} {\sl If $f(t) > 0$ is differentiable, non-increasing, and integrable on $[a,\infty)$, then for almost all $(x,y) \in S \times S$, \[\frac{1}{T}\int_a^T |\tilde{N}(t,x,y)|\,dt \ll e^{T}(Tf(T))^{-1/2}\,.\]For example, taking $f(t)=t^{-1}\log^{-(1+\epsilon)}t$, we obtain, for almost all \mbox{$(x,y)\! \in S\! \times \! S$,} the estimate \[\frac{1}{T}\int_a^T |\tilde{N}(t,x,y)|\,dt \ll e^{T}\log^{\frac{1}{2} +\epsilon}T\,,\]for any $\epsilon > 0$.}

\bigskip

We conclude with a brief account of the situation in general dimension, for which the results in the paper \cite{hill-parnovski} of Hill and Parnovski are very useful. 

As already mentioned, Hill and Parnovski analyze a formulation in which $x$, the center of the ball, varies, but $y$, whose orbit determines the lattice points, is fixed, and they study the variance of $N(T,x,y)$, taken over the $x$ variable. As we have also mentioned, in their version, the variance remains finite in the non-compact finite volume case, since then the integration is taken over $S$, rather than over $S\times S$, and $y$ is thereby constrained from excursions which would otherwise lead to an infinite square integral. Consequently, they discuss the non-compact as well as the compact case, and in the non-compact case, the contribution from the continuous, as well as that from the discrete spectrum must be considered. The Selberg transform plays a role in both cases, which they analyze in terms of the Gaussian hypergeometric function $_2F_1(a,b;c;z)$, an approach which was introduced in this context in a startlingly prescient paper by Delsarte \cite{delsarte} and later used, among others, by Levitan \cite{levitan}. We note that in his 1942 announcement, Delsarte gave what would now be called the Selberg transform for the important case of the indicator function of a ball of radius $T$ in $\mathbb{H}^2$. In that announcement, he also explicitly identified the importance of small eigenvalues for the asymptotics of $N(T,x,y)$, and specifically raised the question of their possible existence. His short {\it Comptes Rendus} announcement \cite{delsarte} does not contain all his detailed derivations of these results, but they are present in in his complete works \cite{delsarte-works}.

We now recall the required results from \cite{hill-parnovski}, which are somewhat intricate to state, and for which we require some preliminary definitions from \cite{hill-parnovski}. Since the notation of this paper differs in a few particulars from that of \cite{hill-parnovski}, we have, in stating what is required, made small necessary adjustments to bring the notation into conformity with that of this paper.

\bigskip

{\sc Definitions:}

\begin{enumerate}

\item $D= (k-1)/2$.

\item $s_j$ is one of the two roots of $\lambda_j = s_j(k-1-s_j)$, all of which are present on the union of the real segment $[0,k-1]$ and the line $\mathrm{Re}(s) = D$. For definiteness, we will take the $s_j$'s to lie on the union of $[D,k-1]$ and the half-line defined by $\mathrm{Re}(s)= D$, $\mathrm{Im}(s)>0$. Note that the definition of the $s_j$'s differs slightly from the Selberg convention for the $r_j$'s, for which the $s_j$'s are counterparts in \cite{hill-parnovski}, and that $r$ and $s$ are related by:\ \ $s = -ir +D$.

\item $w(s) = \pi^D\Gamma(s-D)/\Gamma(s+1)$.

\item \[f_s(x) = w(s)(1-x)^k \sum_{0\leq n < s-D} \frac{(k-s)_n (D+1)_n}{(D-s+1)_n}\,\frac{x^n}{n!}\,.\]

where $(x)_n$ is a Pochhammer symbol, defined, for $n=0$, by $(x)_0 = 1$, and for $n \geq 1$, by $(x)_n = x(x+1)\cdots (x+n-1)$.

\item \[ c_{\log}(s) = \frac{2(-1)^{s-D}\pi^{k-1}}{(s-D)\Gamma(s+1)\Gamma(k-s)}\]if $s-D \in \mathbb{Z}$, and is $0$ otherwise.

\item \[h_k = -\sum_{j=1}^{(k-1)/2} \frac{1}{j}\] if $k$ is odd, and is \[ \log 4 - \sum_{j=0}^{(k-2)/2} \frac{2}{2j +1}\] if $k$ is even.

\end{enumerate}

\bigskip

{\sc Theorem A} (Hill and Parnovski. Theorem 1 of \cite{hill-parnovski}).

\medskip

{\sl Assume $S$ compact and $k$-dimensional, with $k\geq 4$. Then

\[\int_{S} |(N(T,x,y)- B_{k}(T)/V|^2\,dx =  \sum_{j:0<\lambda_j<D^2}f_{s_j}(e^{-2T})^2|\varphi_j(y)|^{2}e^{2s_jT}\]

\[+ \sum_{j:0<\lambda_j<D^2}c_{\log}(s_j)|\varphi_j(y)|^{2}Te^{(k-1)T} + \frac{4\pi^{k-1}}{\Gamma(D+1)^2}\sum_{j:\lambda_{j}=D^2}|\varphi_{j}(y)|^{2}(h_k +T)^2 e^{(k-1)T}\]\[ + O\left(e^{(k-1)T}\right)\,.\]If $k$ equals $2$ or $3$, this specializes to\[\int_{S} |(N(T,x,y)- B_{k}(T)/V|^2\,dx =  \sum_{j:0<\lambda_j<D^2}(w(s_j))^2|\varphi_j(y)|^{2}e^{2s_jT}\]

\[ + \frac{4\pi^{k-1}}{\Gamma(D+1)^2}\sum_{j:\lambda_{j}=D^2}|\varphi_{j}(y)|^{2}(h_k +T)^2 e^{(k-1)T} + O\left(e^{(k-1)T}\right)\,.\]In all cases, the implied constant in the $O$-term depends on $\Gamma$.}

\pagebreak

{\sc Theorem B} (Hill and Parnovski. Theorem 2 of \cite{hill-parnovski}).

\medskip

{\sl Denote, depending on the dimension, the expression on the right side of the last equality, except for the $O\left(e^{(k-1)T}\right)$ term, by $H(T,y)$. Assume $S$ non-compact, of finite volume, and $k$-dimensional. Then, for fixed $y$, \[\int_{S} |(N(T,x,y)- B_{k}(T)/V|^2\,dx = H(T,y) + \frac{\pi^{k-1}}{\Gamma(D+1)^2}|\mathcal{E}(y,D)|^2 T e^{(k-1)T}\]\[+ O\left(e^{(k-1)T}\right)\,,\]where $|\mathcal{E}(y,D)|$ is the Euclidean norm of the vector of the normalized Eisenstein series corresponding to the cusps of $S$, evaluated at $(y,D)$. The implied constant in the $O$-term depends on $y$ and $\Gamma$.}

\medskip

In view of Theorems A and B, bearing in mind that they describe asymptotic behavior for large $T$, it is natural, as before, to subtract the corresponding principal asymptotics from the Fourier expansion of $N(T,x,y)$, and to introduce a function $\tilde{N}(T,x,y)$, defined separately for the compact and non-compact cases, to asymptotically represent the Fourier expansion of the remainder when $T$ is large. Depending on circumstances, the definitions will be valid for $L^2(S)$ or $L^2(S\times S)$.

Beginning with the definition for the compact case, for which we will initially take the context to be $L^2(S\times S)$, and assuming $k\geq 4$, set, for large $T$,

\[\tilde{N}(T,x,y) = N(T,x,y) - B_k(T)/V\] \[- \sum_{j:0<\lambda_j<D^2}\left[f_{s_j}(e^{-2T})e^{2s_jT} +c_{\log}(s_j)Te^{((k-1)T)/2})\right]^{1/2}\varphi_j(x)\varphi_j(y)\]
\[ - \;\;\frac{2\pi^{(k-1)/2}}{\Gamma(D+1)}\sum_{j:\lambda_{j}=D^2}(h_k +T) e^{((k-1)T)/2}\varphi_j(x)\varphi_{j}(y)\,.\]

{\sc Remark.}For $s_j$ in the applicable range and large $T$, the quantity in brackets above is positive, and the positive square root is taken in the definition.

\medskip

If $k$ is $2$ or $3$, set\[\tilde{N}(T,x,y) = N(T,x,y) - B_k(T)/V\]\[-\sum_{j:0<\lambda_j<D^2}w(s_j)e^{s_jT}\varphi_j(x)\varphi_j(y)\]
\[ - \;\;\frac{2\pi^{(k-1)/2}}{\Gamma(D+1)}\sum_{j:\lambda_{j}=D^2}(h_k +T) e^{((k-1)T)/2}\varphi_j(x)\varphi_{j}(y)\,.\]

To facilitate the definition of $\tilde{N}(T,x,y)$ in the non-compact case, we will in all cases denote by $H(T,x,y)$ the quantity subtracted from $N(T,x,y)$ on the right sides of the equations defining $\tilde{N}(T,x,y)$, bearing in mind that $y$ is held fixed in the non-compact case, and that for fixed $T$ and $y$, $N(T,x,y)$ is piecewise constant and of compact support on $S$. 

Specifically, in the non-compact case, we set \[\tilde{N}(T,x,y) = N(T,x,y) - H(T,x,y) - C(T,x,y)\,,\]where the new  term $C(T,x,y)$ arises from the role played by Eisenstein series in the spectral analysis of $S$ for this case. For the purposes of this paper, we will only define $C(T,x,y)$, but it originates from considerations which arise from Theorem B and the literature on Eisenstein series, to which we refer the reader seeking further details (cf.,\ for example, \cite{iwaniec1}, \cite{mueller}). 

In order to define $C(T,x,y)$, we recall that for suitable functions $f(x)$ on $S$, the Fourier-Eisenstein transform corresponding to a cusp is defined by\[ \hat{f}(D+it) = \int_S f(x)\overline{\mathcal{E}}_j(x,D+it)\,dx\,,\]where $\mathcal{E}_j(x,s)$ is the normalized Eisenstein series associated to the cusp. The domain of $\hat{f}$ is the line $D+ it$ ($-\infty < t < \infty$) in the complex plane, which, in the previously described parametrization by $s$ of eigenvalues, corresponds to $\lambda \ge D$. By the general theory, the continuous spectrum associated to such a cusp lies on this line, and the projection of $f$ on the orthocomplement in $L^2(S)$ of the subspace spanned by eigenfunctions corresponding to the discrete spectrum is given by the sum of the inverse transforms\[\frac{1}{4\pi}\sum_j\int_{-\infty}^{\infty} \hat{f}(D+it)\mathcal{E}_j(x,D+it)\,dt\,.\]Now, regarding $y$ and $T$ as constants, so that $\hat{N}$ can be considered as a function of $D+it$, we define\[C(T,x,y) = \frac{1}{4\pi}\sum_j\int_{-\infty}^{\infty} \hat{N}(D+it)\mathcal{E}_j(x,D+it)\,dt\,.\]

With the preceding definitions, examination of the square of the $L^2$ norm of \linebreak $\tilde{N}(T,x,y)$ now brings Theorem B into play. We remark that the Eisenstein term in Theorem B derives from an analysis of a dominating contribution, for fixed $y$, to the $T$ asymptotics of the square of the $L^2$ norm of $\hat{N}$ (cf.\ Theorem 4 of \cite{hill-parnovski})).

This said, it now follows from Theorems A and B (expand the integrand on the left in the statements of those theorems) that in both the compact and non-compact cases,\[\int_S |\tilde{N}(T,x,y)|^2\,dx = O\left(e^{(k-1)T}\right)\,,\]with the previously described dependencies of the implied constants.

Note, as is mentioned in \cite{hill-parnovski}, that in the compact case, in which we can simultaneously vary $x$ and $y$, we can further integrate this with respect to $y$, to deduce that in the compact case, we also have\[\int_{S\times S} |\tilde{N}(T,x,y)|^2\,dxdy = O\left(e^{(k-1)T}\right)\,.\]

With these preparatory facts in hand, we can now state various Hardy-type integral estimates for the hyperbolic case. Their derivation exactly follows the already described method, so we will content ourselves with statements of the conclusions.

\medskip

{\sc Theorem 3.} {\sl Suppose that $S$ is compact, and that $f(t) > 0$ is differentiable, non-increasing, and integrable on $[a,\infty)$. Then for almost all $(x,y) \in S \times S$, \[\frac{1}{T}\int_a^T |\tilde{N}(t,x,y)|\,dt \ll e^{\frac{1}{2}(k-1)T}(Tf(T))^{-1/2}\,.\]For example, taking $f(t)=t^{-1}\log^{-(1+\epsilon)}t$, we obtain, for almost all \mbox{$(x,y)\! \in S\! \times \! S$,} the estimate \[\frac{1}{T}\int_a^T |\tilde{N}(t,x,y)|\,dt \ll e^{\frac{1}{2}(k-1)T}\log^{\frac{1}{2} +\epsilon}T\,,\]for any $\epsilon > 0$.}

\medskip

{\sc Theorem 4.} {\sl Suppose that $S$ is either compact, or non-compact with finite volume, and that $f(t) > 0$ is differentiable, non-increasing, and integrable on $[a,\infty)$. Then, assuming $y$ fixed, for almost all $x \in S$, \[\frac{1}{T}\int_a^T |\tilde{N}(t,x,y)|\,dt \ll e^{\frac{1}{2}(k-1)T}(Tf(T))^{-1/2}\,.\]For example, taking $f(t)=t^{-1}\log^{-(1+\epsilon)}t$, we obtain, for almost all $x \in S$, the estimate \[\frac{1}{T}\int_a^T |\tilde{N}(t,x,y)|\,dt \ll e^{\frac{1}{2}(k-1)T}\log^{\frac{1}{2} +\epsilon}T\,,\]for any $\epsilon > 0$.}

\medskip

{Remark. }For $k=2$, this accords well, up to the logarithmic factor introduced by our methodology, with the numerical results described in \cite{phillips-rudnick}.

\newcommand{\noopsort}[1]{}

\bigskip

\bigskip

{\sc \noindent Ph.D.\ Program in Mathematics\\CUNY Graduate Center\\365 Fifth Avenue, New York, NY 10016}

\end{document}